\documentclass{amsart}
\usepackage{graphicx}
\usepackage{amsmath}

\newtheorem{thm}{Theorem}

\newtheorem{prop}{Proposition}
\newtheorem{lem}{Lemma}

\newtheorem{defn}{Definition}

\newcommand{\Ker}{\mbox{Ker}}

\begin{document}

\title{The Root Closure of a Ring of Mixed Characteristic}

\author{Paul C. Roberts}
\email{roberts@math.utah.edu}
\address{Department of Mathematics, University of Utah, 155 S 1400 E, Salt Lake city, Utah 84112-0090}
\thanks{The author was supported by NSF grant 0500588.}

\begin{abstract}
We define a closure operation for rings of mixed characteristic and verify that the
closure is a ring. We then show that this closure produces a ring with
good properties with respect to its Fontaine ring and give an example to
show that rings that are not closed in this sense do not satisfy
these properties.
\end{abstract}

\maketitle

\section{Introduction}
The results in this paper are part of a program  to understand
rings of mixed characteristic by studying 
their associated Fontaine rings. We will define Fontaine rings
and outline the main properties that we use in Section 3; we mention here that
they give a ring of positive characteristic from which, under
certain conditions, the original ring can be reconstructed up to $p$-adic
completion.  The difficulty comes from the ``certain conditions"  that have to
apply;  essentially what is necessary is that there are enough $p$th roots
in the ring.  For absolutely integrally closed rings, for example,  this works well, but
for Noetherian rings it does not work at all, and adjoining $p$th roots of all elements
produces  a huge extension that is difficult to deal with.  In this paper
we describe a much smaller extension, which we call the root closure,
that makes the reconstruction of the $p$-adic completion from the
Fontaine ring work correctly.

In section 2 we give the basic definitions and elementary properties
of the  root closure.  In section 3 we desribe the connection
with Fontaine rings, and in section 4 we give an example to show
how this works in practice.

\section{Basic Definitions}  Let $p$ be a prime number, and let $R$ be a 
commutative ring (with 1) of mixed characteristic $p$.
The only assumption we make is  that $p$ is not a zero divisor in $R$;
in most cases of interest $R$ is a quasi-local domain
and $p$ is a nonzero element of its maximal ideal.
However, we do not even exclude the case where $p$ is
a unit, although the ring does not actually have mixed characteristic
and the construction  is not interesting in that case.  
We do not assume that $R$ is Noetherian.

Let $R_p$ denote the localization of $R$ obtained by inverting $p$; our
assumption implies that $R\subseteq R_p$.  

\begin{defn} The root closure of $R$, denoted $C(R)$, is the set
of all $x\in R_p$ such that $x^{p^n}$ is in $R$ for some integer $n\ge 0$.
\end{defn}

We note that if $x^{p^n}\in R$, then $x^{p^m}\in R$ for
all $m\ge n$, and in fact we also have that $x^{kp^n}\in R$ for
all positive integers $k$.  

\begin{prop}
If $R$ is a ring such that $p$ is not a zero divisor on $R$, then $C(R)$ is a subring of $R_p$.
\end{prop}

{\bf Proof.}  We note first  that if we have a finite number of elements $s_i$
in $C(R)$ we can assume that there a common $n$ such that  $s_i^{p^n}$
is in $R$ for all $i$.   Since $C(R)$ is by definition contained
in $R_p$, we can also assume that there is a common integer
$k$ such that $p^ks_i\in R$ for all  $i$.

It is clear that $R\subseteq C(R)$ (so in particular $1\in C(R)$), and that if $s$ and $t$ are in
$C(R)$, with $s^{p^n}$ and $t^{p^n}$ in $R$, then $(st)^{p^n}=s^{p^n}t^{p^n}\in R$, so
$st\in C(R)$.
It is also clear that if $s\in C(R)$, then $-s\in C(R)$.
 
We must now show that if $s$ and $t$ are in $C(R)$, 
then $(s+t)^{p^N}\in R$ for some integer $N$.  Let $n$  be a positive integer
such that $s^{p^n}$ and $t^{p^n}$ are in $R$ and let $k$ be a positive integer such
that $p^ks$ and $p^kt$ are in $R$.

We first prove an elementary lemma on binomial coefficients.  We use the notation $p^r\parallel m$
to mean that $p^r$ is the highest power of $p$ that divides $m$.

\begin{lem} Let $m$ be a positive integer, and let $i$ be an integer with $1\le i\le p^m$.
If $p^r\parallel i$, then $p^{m-r}\parallel  {{p^m}\choose i}$.
\end{lem}

{\bf Proof.}  We use induction on $i$.  
 For $i=1$, ${{p^m}\choose i}=p^m$, $r=0$, and $p^{m-r}=p^m\parallel
 p^m$, so the result is correct.  
 
 We now assume that
$2\le i\le p^m$ and assume the result is true for $i-1$.  Since
$${p^m\choose i}={p^m\choose {i-1}}\left(\frac{p^m-i+1}{i}\right),$$
the only time the power of $p$ that divides $i$ or ${p^m\choose i}$
will change
is if $p|i$ or $p|(p^m-i+1)$, which means that $p|i-1$.  

If $p|i$, then $p\not | i-1$,
so the induction hypothesis implies that
$p^m\parallel {{p^m}\choose {i-1}}$.  Then if $p^r\parallel i$, 
 to obtain $ {p^m\choose i}$ from $ {p^m\choose i-1}$
we multiply by a number prime to $p$ and divide by $i$, so we
conclude that $p^{m-r}\parallel p$.

If $p|i-1$, and if $p^r\parallel i-1$, then
by induction $p^{m-r}\parallel
{ {p^m\choose i-1}}$, and a similar computation shows that 
$p^m\parallel {p^m\choose i}$; since $p\not| i$
in this case, this proves the result.  \hskip1in$\Box$

\bigskip

We now return to the proof of the theorem.  As above, let $s^{p^n}$,  $t^{p^n}$,
 $p^ks$, and $p^kt$ be in $R$.  We claim
that if $N>2kp^n+n$, then $(s+t)^{p^N}\in R$.

We have 
$$(s+t)^{p^N} = \sum_{i=0}^{p^N}{{p^N}\choose i}s^it^{p^N-i}.$$

We claim that every term in the sum on the left hand side is in $R$.  First, if $p^n$ divides $i$,
then $p^n$ divides $p^N-i$, and both $s^i$ and $t^{p^N-i}$ are in $R$.
On the other hand, if $p^n$ does not divide $i$, then the above lemma implies
that $p^{N-n}$ divides ${{p^N}\choose i}$.  By our choice of $N$, this implies
that $p^{2kp^n}$ divides ${p^N}\choose i$.  

Write $i=ap^n+u$ and $p^N-i=bp^n+v$, with $u$ and $v$ integers such that $1\le u,v< p^n$.
We have $s^i=s^{ap^n+u}=s^{ap^n}s^u.$  As noted above, $s^{ap^n}\in R$.
Since $p^ks\in R$, we have $p^{ku}s^u\in R$, so, since $u<p^n$, $p^{kp^n}s^u\in R$.
Thus $p^{kp^n}s^i\in R$.  Similarly, $p^{kp^n}t^{p^N-i}\in R$.  Thus
$p^{2kp^n}s^it^{p^N-i}\in R$, so, since $p^{2kp^n}$ divides ${p^N}\choose i$,
${{p^N}\choose i}s^it^{p^N-i}\in R$.  Thus we have shown that every term
in the sum is in $R$, so $(s+t)^{p^N}\in R$ and $s+t\in C(R)$.  Thus $C(R)$ is
a ring.

\bigskip

If $R=C(R)$, we say that $R$ is  root closed.  We note that this is a weaker
statement than saying that it is closed under taking $p$th roots; we also note that if
$S$ is an arbitrary extension of $R$, it is not true that the set of elements $x$ of $S$
such that $x^{p^n}$ is in $R$ for some $n$ forms a ring.

\section{The  root closure and Fontaine rings.}

Let $R$ be a ring of mixed characteristic as above. We define
the {\em Fontaine ring} of $R$, which we denote $E(R)$, by
$$E(R)=\lim_{\leftarrow} R_n,$$
where each $R_n$, defined for integers $n\ge 0$, is $R/pR$, and the
map from $R_{n+1}$ to $R_n$ is the Frobenius map.
An element of this ring is thus given by a sequence $r_0,r_1,\ldots $ of
elements of $R/pR$ with $r_{n+1}^p=r_n$ for all $n\ge 1$.  We
denote this sequence $(r_n)$.  

It is rather clear from the definition that if there are not very many elements that have $p^n$th roots
modulo $p$, the Fontaine ring will be very small.  However, the
only assumption that we make on $R$ is that it contain a $p^n$th
root of $p$ for each $n$.  In this case a compatible specific 
choice of $p^n$th root for each $n$ defines an element $(p^{1/p^n})$
of $E(R)$;
we denote this element $P$.

Fontaine rings have been studied for valuation rings in connection with Galois representations 
by Fontaine \cite{F}, Wintenberger\cite{W}, and others. They
have been studied for more general rings by Andreatta \cite{A}.  
For their basic properties we refer to
Gabber and Ramero \cite{GR} Section 5.2; our notation is taken
essentially from that source.   We note that 
there are alternative definitions for this ring and that there are other
``Fontaine rings" rings defined from this one. 

We recall (see \cite{GR}) that $E(R)$ is a perfect ring of
characteristic $p$; in fact, the $p$th root of $(r_n)$
is simply $(s_n)$, where $s_n=r_{n+1}$ for
each $n$. 

One of the most useful properties of Fontaine rings is that the
$p$-adic completion of $R$ can be reconstructed from $E(R)$
and, for certain rings, this can be done in a simple way.  More
precisely, there is a map from the ring of Witt vectors on $E(R)$
to the $p$-adic completion of $R$.  We recall that the
$p$-adic completion of $R$ is $\displaystyle\lim_{\leftarrow}R/p^nR$, which
we denote $\hat R$.  We will show that to determine the
kernel of this map it suffices
that $R$ be  root closed.

  Let $W(E(R))$ denote the ring of Witt
vectors on $E(R)$ and let $\hat R$ denote the $p$-adic
completion of $R$.  We refer to Bourbaki \cite{BB} for
general properties of Witt vectors.

We have a map $\tau_R$
from $E(R)$ to $W(E(R))$; it sends $a$ to $(a,0,0,\ldots)$.
The  map $\tau_R$ preserves multiplication but
not addition.  There is a map $u_R$ from $W(E(R))$ to $\hat R$
such that for an element $(r_n)$ of $E(R)$ we have
$$u_R(\tau_R(r_n)) = \lim_{n\to \infty}r_n^{p^n}.$$
Finally, $u_R$ induces a ring homomorphism $\overline u_R$ from
$E(R)/PE(R)$ to $\hat R/p\hat R = R/pR$ that coincides
with the map defined by sending $(r_n)$ to $r_0$.

\begin{lem}\label{kernel1}  Let $F$ be the Frobenius map 
on $R/pR$, where $R$ is as above.  If $R$ is  root closed, then the kernel of $F^n$
is generated by $p^{1/p^n}$ for all $n$.
\end{lem}
Suppose that $a\in R$ and that $a^{p^n}$ is a multiple of $p$, so
$a^{p^n}=pb$ for some $b\in R$.  
Now we have $(a/p^{1/p^n})^{p^n}=pb/p=b\in R$, so
$a/p^{1/p^n}$ is in $C(R)$.  Since $R$ is  root closed, $a/p^{1/p^n}$
is in $R$, so $a\in p^{1/p^n}R$, as was to be shown.

\bigskip

We next show that if $R$ is  root closed, we can determine the kernel of 
the map $\overline u_R$ from $E(R)$ to $R/pR$. 
\begin{prop}
Suppose that $R$ is root closed.  Then $\overline u_R:
E(R)/PE(R)\to \hat R/p\hat R = R/pR$ is injective.
\end{prop}
{\bf Proof.}  Let $R=(r_n)$ be an element of $E(R)$ that goes to zero
in $R/pR$, which means that $r_0=0$ in $R/pR$.
Since $r_n^{p^n} = r_0$, the above lemma implies that  $r_n\in p^{1/p^n}R$ for all $n$.
Let $r_n=p^{1/p^n}s_n$ for each $n$.  

It is clear that letting $S=(s_n)$ would give $R=PS$, but we do not
know that  $S$ is an element of $E(R)$;
that is, that $s_{n+1}^p=s_n$ in $R/pR$ for all $n\ge 1$.  We lift the $s_n$ to
elements of $R$, also denoted $s_n$. In our notation we will
henceforth use congruence modulo $pR$ to denote equality in $R/pR$. 

Although we do not know that
$s_{n+1}^p\equiv s_n$ modulo $pR$,  we do know
that $(p^{1/p^{n+1}}s_{n+1})^p\equiv r_{n+1}^p\equiv r_n \equiv 
p^{1/p^n}s_n$ modulo $p$, which
implies that $s_{n+1}^p\equiv s_n$ modulo $p^{1-1/{p^n}}R$
(since $p$ is not a zero-divisor on $R$).  
We claim that if we let $t_n=(s_{n+1})^p$
for each $n$, then we will have $t_{n}^p\equiv t_{n-1}$ modulo
$p$ for each $n\ge 1$ and $(r_n)=P(t_n)$.  

Since $s_{n+1}^p\equiv s_n$ modulo $p^{1-1/{p^n}}R$, there
is an element $v_n$ of $R$ with 
$$s_{n+1}^p =  s_n + p^{1-1/{p^n}}v_n.$$
If we raise this equation to the $p$th power we obtain
$$s_{n+1}^{p^2} \equiv  s_n^p + (p^{1-1/{p^n}}v_n)^p \equiv  s_n^p$$
modulo $pR$.  This says that $t_{n}^p\equiv t_{n-1}$ modulo $pR$, as
was to be shown.

It remains to show that $r_n\equiv p^{1/p^n}t_n$ modulo $pR$ for each $n$. In fact, we have
$$r_n \equiv r_{n+1}^p\equiv (s_{n+1}p^{1/p^{n+1}})^p\equiv t_np^{1/p^n}\;\;\mbox{modulo}\;\;pR.$$ 

Thus $(r_n)\in PE(R).$
\bigskip

We now come to the main point, which is to describe the 
kernel of the map $u_R$ from
$W(E(R))$ to $\hat R$.

\begin{thm}  Suppose $R$ is root closed.  Then the kernel of
$u_R$  is generated by $P-p$.
\end{thm}

{\bf Proof.} We have a diagram

$$\begin{array}{ccccccccc}
0&\to&\Ker(u_R)&\to &W(E(R))&\stackrel{u_R}{\to} &\hat R\\
&&\downarrow&&\downarrow&&\downarrow\\
0&\to&\Ker(\overline u_R)&\to &E(R)&\stackrel{\overline u_R}{\to} &R/pR
\end{array} $$

The map from $W(E(R))$ to $E(R)$ is reduction modulo $p$.

Let $x$ be an element in the kernel of  $u_R$.  Mapping $x$ down to $E(R)$
we get an element of the kernel of $\overline u_R$, so by Lemma
\ref{kernel1} we obtain a multiple of $P$, which we write $eP$.  Lift $e$ to an element
of $W(E(R))$, say $w_1$; we then have that $x-Pw$ goes to
zero in $E(R)$, so $x-Pw_1\in pW(E(R)).$  Thus $x-(P-p)w_1$  
is in $pW(R)$, so we can write 
$$x = (P-p)w_1 + pv_1    \eqno{(*)}$$
for some $v_1\in W(E(R)).$

To complete the proof we use induction on $k$.

Suppose we have $k\ge 1$ and $w_k$ and $v_k$ in $W(E(R))$ with
$$x=(P-p)w_k+p^kv_k.$$  Equation ($\ast$) gives the case $k=1$.

Since $x$ and $P-p$ are in the kernel of $u_R$, we have
$$p^ku_R(v_k)=u_R(p^kv_k)=u_R(x-(P-p)w_k)=0.$$  Thus, since
$p$ is a nonzerodivisor in $\hat R$, $v_k$ is in the kernel of $u_R$.   Using the
same argument that we used above for $x$, 
we can write

$$v_k=(P-p)y+pv_{k+1}$$
for some $y$ and $v_{k+1}$ in $W(E(R))$.  We then have
$$ x=(P-p)w_k+p^kv_k=(P-p)w_k+p^k((P-p)y+pv_{k+1})$$
$$=(P-p)(w_k+p^ky)+p^{k+1}v_{k+1}.$$
Letting $w_{k+1}=w_k+p^ky$ we have 
$$x=(P-p)w_{k+1}+p^{k+1}v_{k+1}$$  
with $w_{k+1}\equiv w_k$ modulo $p^k$.  Since $W(E(R))$ is complete in the $p$-adic
topology, if we let $w$ be the limit of the $w_k$, we have
$x=(P-p)w$.  Thus the kernel of $u_R$ is generated by $P-p$.

\section{An Example}

We give a simple example to illustrate the constructions described above.
Let $V_0$ be the ring of $p$-adic integers for some prime $p>3$, and let
$R_0=V_0[[x,y]]/(p^3+x^3+y^3).$  Let $R$ be the ring obtained by
adjoining $p^n$th roots of $p,x$, and $y$; specifically, we adjoin elements
$\pi_n,x_n,$ and $y_n$ with $\pi_0=p,x_0=x,$ and $y_0=y$ and such that
$\pi_{n+1}^p=\pi_n$ and similarly for the $x_n$ and $y_n$ for each
$n\ge 1$. 

To describe the ring $R$, we first define a ring $S$ similarly,
letting $S_0$ be the power series ring $V_0[[x',y']]$
and adjoining $p^n$th roots $\pi'_n,x'_n,y'_n$ following
the same procedure as for $R$.  We let $S$ be the union of
the $S_n=S_0[\pi'_n,x'_n,y'_n]$; in this case each $S_n$ is a regular local ring.
There is a map from $S$ to $R$ that sends $\pi'_n$ to $\pi_n$,
$x'_n$ to $x_n$, and $y'_n$ to $y_n$.

We claim that the kernel of this map is generated by $p^3+x'^3+y'^3$.
To show this it suffices to show that $p^3+x'^3+y'^3$ is prime in
$S_n$ for each $n$.  In $S_n$ this polynomial can be written
as a polynomial in $y'_n$ as
${y'}^{3p^n}+({\pi'}_n^{3p^n}+{x'_n}^{3p^n})$, which is prime since
${\pi'}_n^{3p^n}+{x'}_n^{3p^n}$ is a product of distinct prime elements
of $S_n$ (using Eisenstein's criterion, for example).

 Let $P,X,$ and $Y$ be the elements $(\pi_n),(x_n)$, and
$(y_n)$ of $E(R)$. Consider the element $\eta=(r_n)=P^3+X^3+Y^3$.  Its
zeroth component $r_0$ is $p^3+x^3+y^3=0$, so $\eta$ is in the
kernel of $u_R$.  We claim that $\eta$ is not in $PE(R)$.
If it were, its $r_1$ component, $p^{3/p}+x^{3/p}+y^{3/p}$
would have to be in the ideal generated by $p^{1/p}$,
which means that the corresponding power series 
$p^{3/p}+{x'}^{3/p}+{y'}^{3/p}$
would be in the ideal generated by $p^{1/p}$ and 
$p^3+x'^3+y'^3$.  This is clearly not the case,
so $\eta\not\in PE(R)$.

On the other hand, it is easy to see that the
elements $(p^{3/p^n}+x^{3p^n}+y^{3p^n})/p^{1/p^n}$ are in
$C(R)$ and that $\eta\in PE(C(R))$.

We remark that it can be shown that $E(R)$ in this example
is a completion of a power series ring over a field in
three variables, so that an attempt to recover $\hat R$
by taking $W(E(R))/(P-p)$ would give a ring of dimension
3 rather than 2, the dimension of $R$.

\end{document}